# Recurrent Theme of Pick's Theorem

Jacek M. Kowalski

(University of North Texas, Denton *)

ABSTRACT. We review and possibly add some new variant to the existing derivations of the formula for the area of Jordan lattice polygons drawn on two-dimensional lattices. The formula is known as Pick's theorem and is related to the number theory elementary result- Bézout lemma. It is pointed out that Euclidean algorithm can be easily used in construction of infinite number of distinct primitive cells for any two-dimensional lattice. Pick formula itself can also be obtained in an elementary "cut and re-assemble" finite process.

## 1.Introduction

Every physicist's square lattice may be identified with 2-D integer lattice $\mathbf{Z}^2$, provided that the lattice constant is set to 1. *Lattice polygons* can be drawn on such lattice by joining lattice points. In 1899 Einstein's friend Georg Pick, working on connection between geometry of arbitrary two-dimensional lattices and number theory, published an important paper [1], leading, among other results, to an expression for the area $A(P)$ of a square lattice polygon $P$:

$$A(P) = \frac{1}{2} N^b(P) + N^i(P) - 1, \tag{1}$$

where $N^b(P)$ is the number of lattice points on the boundary of $P$ and $N^i(P) \geq 0$ is the number of "internal" lattice points enclosed by $P$. This result is typically quoted as Pick's theorem.

An immediate conclusion of the Pick's formula (1) is the fact that any parallelogram formed by square lattice vectors such that: i) its boundary contains only four lattice points ($N^b$=4) and ii) does not enclose internal lattice points ($N^i$=0), must be a parallelogram of unit area. Additionally, any such parallelogram can be used to generate the entire square lattice.

The latter fact was proved in an independent, elegant, but taking almost a page, proof in the classical text *Geometry and the Imagination* [2], for any lattice parallelogram which "does not have any lattice points in its interior or on its boundary other than its vertices". Sources [3] indicate that Hilbert got to know Pick personally, but Pick's formula hasn't been used or quoted in [2]. (As a side remark, let's also note that Hilbert & Cohn-Vossen's result could still be quoted in quite imprecise way, as in [4]: "Any parallelogram on the lattice which two opposite sides have length 1 has unit area").

There are many proofs of the Pick's result published both in books, journals and on independent web pages [5-20]. Though many authors duly quote, difficult to access for a long time, Pick's original paper published in an obsolete journal, but no one of them (among known to the author of this note) discusses main line of Pick's original reasoning. In fact, the "area formula" was used in [1] only as a tool in proving a statement from elementary number theory, known nowadays as Bezout's identity.

In section 2, after some necessary elementary preliminaries, we review variants of the "Pick's theorem" proof, and advertise one more which explicitly uses a special case of Bezout's identity as a starting point. The related well-known Euclid algorithm allows effective construction of lattice unit parallelograms of arbitrary extension. As already noticed by many authors, many of the steps in the proof can be successfully presented to, or formulated as exercises for, students with minimum preparedness, but with some sensitivity to the beauty of mathematical reasoning. This is also a secondary aim of this overview, with less formal presentation style and avoidance of some technicalities.

In Section 3 we briefly discuss Pick's original approach which, in essence, was as a reverse path, "from geometry to an elementary number theory".

The closing section 4 contains some further examples and suggestions for presenting the topic at, sometimes, even more elementary level, and comments on relations of the Pick's Theorem with more general results like surveyor's formula for polygons. Finally, one can also find there some references to quite interesting, more recent generalizations of Pick's formula.

…

## 2. **Preliminaries and variants of the proof**.

a) Two -dimensional lattices and their primitive cells

Let $\{\mathbf{p},\mathbf{q}\}$ be a basis in $\mathbb{R}^2$. A lattice $\Lambda(\mathbf{p},\mathbf{q})$ is a subgroup of $\mathbb{R}^2$ of all linear combinations of basal vectors with integer coefficients:

$$\Lambda(\mathbf{p},\mathbf{q}) = \{m\boldsymbol{p} + n\boldsymbol{q} \mid (m,n) \in \mathbb{Z}^2\} \quad . \qquad (2)$$

Every lattice generates a tiling of $\mathbb{R}^2$ by parallelograms of area $A(\Lambda) = |\mathbf{p}\times\mathbf{q}|$ (called primitive or fundamental cells).Lattices constructed from different bases $\{\mathbf{p},\mathbf{q}\}$ and $\{\mathbf{r}, \mathbf{s}\}$ may coincide as sets of points in $\mathbb{R}^2$. This happens if and only if two bases are related by a linear transformation $A$ with integer matrix elements $a, b, c, d$ and such that $|\det A| = 1$:

$$\begin{aligned} \mathbf{r} &= a\mathbf{p} + b\mathbf{q}, \\ \mathbf{s} &= c\mathbf{p} + d\mathbf{q}. \end{aligned} \qquad (3)$$

It follows from (3) that primitive cells formed by $\mathbf{p},\mathbf{q}$ and $\mathbf{r}, \mathbf{s}$ vectors have the same area:

$$A(\Lambda) = |\mathbf{r} \times \mathbf{s}| = |\mathbf{p} \times \mathbf{q}| \quad . \qquad (4).$$

Area of a primitive cell can be selected as a unit area of a given tiling or as equal two in units of areas of two congruent triangles forming the primitive cell (as it was done in [1]).

As long as special point symmetries of a given lattice are of no concern, all 2-D lattices are isomorphic to $\mathbb{Z}^2$, which in turn can be considered as square lattice spanned on ***i***=(1,0) and ***j***=(0,1) with primitive squares of area 1. In this base all pairs of vectors with components (*a,b*) and (*c,d*) satisfying the condition: $|ad\text{-}cb| = 1$ will determine a primitive cell of area 1. This observation, as we will see below, allows to construct infinitely many distinct primitive parallelograms with vertices at lattice points all having "unit area". To show that, and comment on Pick's and related results, we should introduce the related terminology.

. We will call a lattice vector *simple* if representing it oriented segment does not pass through any additional lattice points. It is easy to check that a lattice vector is simple if and only its components are coprime. A *chain of lattice vectors* is a sequence of simple lattice vectors where the ending point of the predecessor is also the origin of the subsequent vector. A *simple path* on square lattice is a chain of non-overlapping simple lattice vectors where all ending points are visited just once. A *simple*, o*riented lattice polygon* is a closed simple path of non-intersecting simple lattice vectors not sharing their internal points (which can also be called an *oriented Jordan polygon*).

Below we will consider only simple lattice polygons in the above sense. Number of lattice points on a simple polygon *P* will be denoted, as above, by $N^b(P)$ -it clearly coincides with the number of simple lattice vectors in the related chain. All lattice points of *P* where two subsequent simple lattice vectors change their spatial orientation are *vertices* of P. Number of vertices of *P* does not exceed $N^b(P)$ but it is always larger or equal three. As before, $N^i(P) \geq 0$ is the number of lattice points enclosed by *P* but not belonging to *P*. For any lattice polygon, the affine form *F*(*P*):= $\frac{1}{2}N^b(P) + N^i(P) - 1$ is an integer or half-integer characteristic of the polygon, always greater or equal ½ (as $N_b \geq 3$). This function has the additive property for unions of polygons not sharing internal points and having a shared simple path on the boundary:

$$F(P_1) + F(P_2) = F(P) . \qquad (5)$$

For visualization of this truly "*Anschaulige Geometrie*" fact see Fig.1

The underlying lattice selected there is a square lattice, and used below terms: "rectangle", right triangle" relate to this particular lattice.

It follows from the figure inspection, that boundary points are ultimately always properly counted with weighting factor ½ in the union, but possible internal *additional* lattice points on a shared dividing path become additional internal lattice points. This additivity property is typically illustrated by a special case of rectilinear dividing paths. Let's also note, in passing, that Pick himself proved and used the additivity property for the function 2$F$($P$) (see Section 4). In modern literature, Pick's formula was first popularized by Steinhaus [5], and the additivity property proof can be found in another classical text [6].

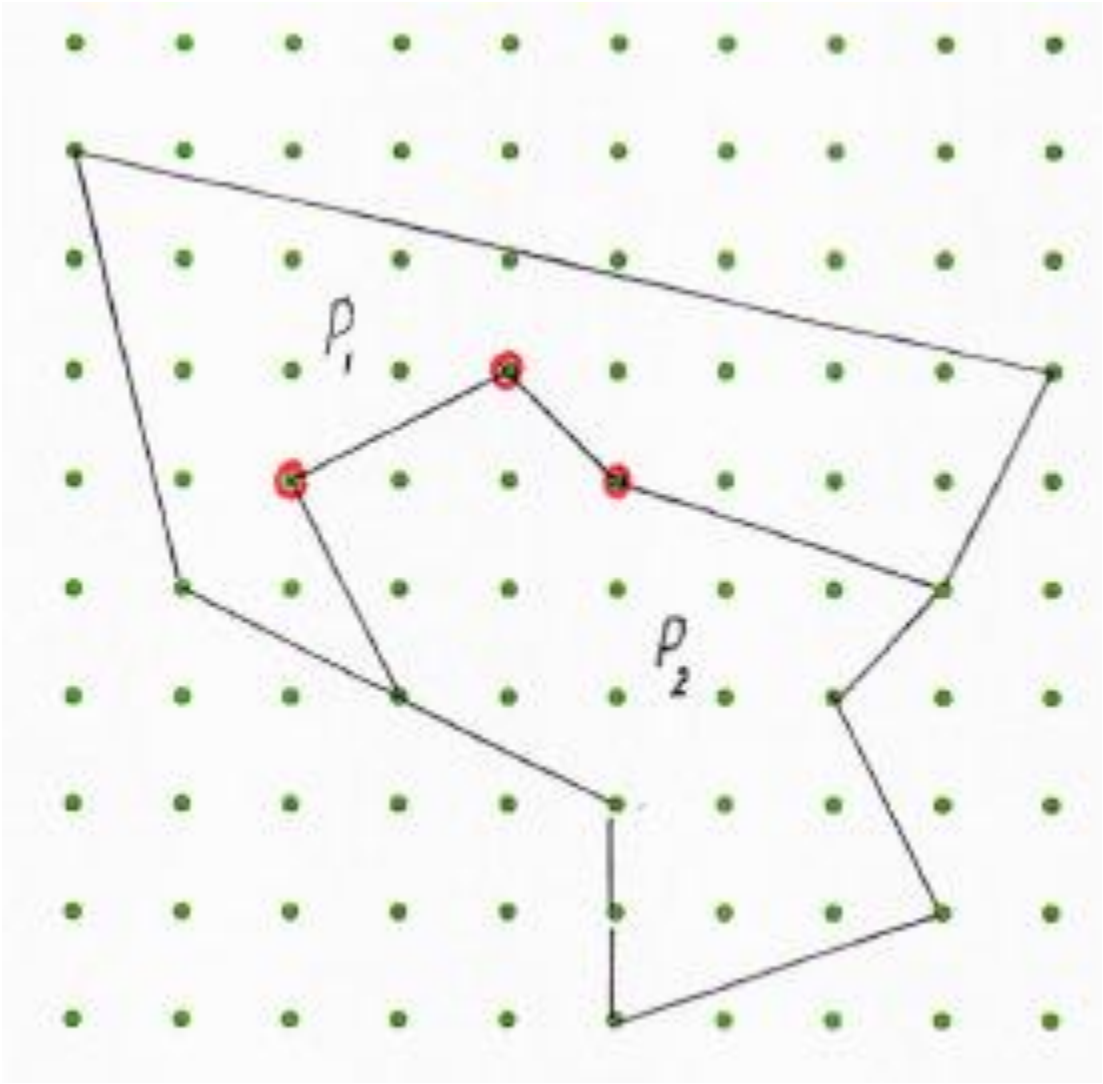

Fig.1

Second, obvious observation is that $F(P)$ is an invariant for all lattice polygons congruent via translations with $P$, as well as for polygons generated by given lattice point symmetry operations.

Further, typical steps in most proofs reduce to the identification of $\frac{1}{2}N^b(P)+N^i(P)-1$ values with areas of respective lattice polygons.

Thus the validity of the Pick's formula is first verified for any lattice rectangle with side lengths $m \geq 1$, $n \geq 1$: $A=mn$, $N^b=2(m+n+1)$, $N^i=(m-1)(n-1)$.

Considering any such rectangle as a union of two congruent right triangles one justifies the theorem for any lattice right triangles having sides along lattice lines (and, as a "side result" not needed in the proof, for unions of such congruent right triangles sharing a side). Finally, any lattice triangle can be inscribed into a lattice rectangle (with all three triangle vertices on the rectangle's perimeter). By "cutting corners" of such rectangle, one validates the theorem for arbitrary lattice triangles (see Section 4a.). Consequently, the theorem is satisfied for any lattice polygon which can be triangulated into or, equivalently, built from lattice triangles.

Another approach to allowing to identify $A(P)$ with the area of the related polygon may consist, first, in a triangulation of a given lattice Jordan polygon into a set of lattice triangles. This can be easily done for a convex polygon or, more generally, by noting that any Jordan n-gon in two-dimensional space can be triangulated into n-2 triangles (see, e.g., a nice online text [19]). One then continues sub- triangulation of triangles in the partition, if they still contain at least one internal point and/ or additional lattice points on their boundary, into three or two "smaller" triangles. This process is continued, if necessary, until the whole polygon is covered by the union of "elementary" triangles with all sides identifiable with simple lattice vectors and no internal points. Clearly, resulting partitions are, in general, not unique, but this fact plays no role in the final steps of the argument. Provided that we can prove that all elementary triangles must have area equal ½, we verify Pick's formula for an elementary triangle and use the additivity property. Let us also note that appealing to the additive property is not necessary in the elementary triangle approach, as the Pick formula can be obtained by a "puzzle- like" process, even more intuitive at the introductory level. Imagine given lattice polygon triangulated, as above, into elementary triangles. Next, start with an arbitrary elementary triangle ( $N^b$=3, $N^i$=0) and keep adding neighboring triangles from remaining ones in the union. At each subsequent steps one either adds one more boundary point with the count of internal points unchanged, or, in some cases, creates an additional internal point by "fitting in" respective elementary triangle into wedge-like regions between two elementary triangles, without increasing the number of boundary points in resulting polygon. For illustration, see Fig. 2. A polygon re-constructed there already contains some less obvious elementary triangles built from indicated vectors: $\triangle$ { (2,5), (1,4)}, $\triangle$ {(3,2),(4,3)} and $\triangle$ {(3,2) and (2,1)}. Green or red triangles are just about to be added at steps of the process.

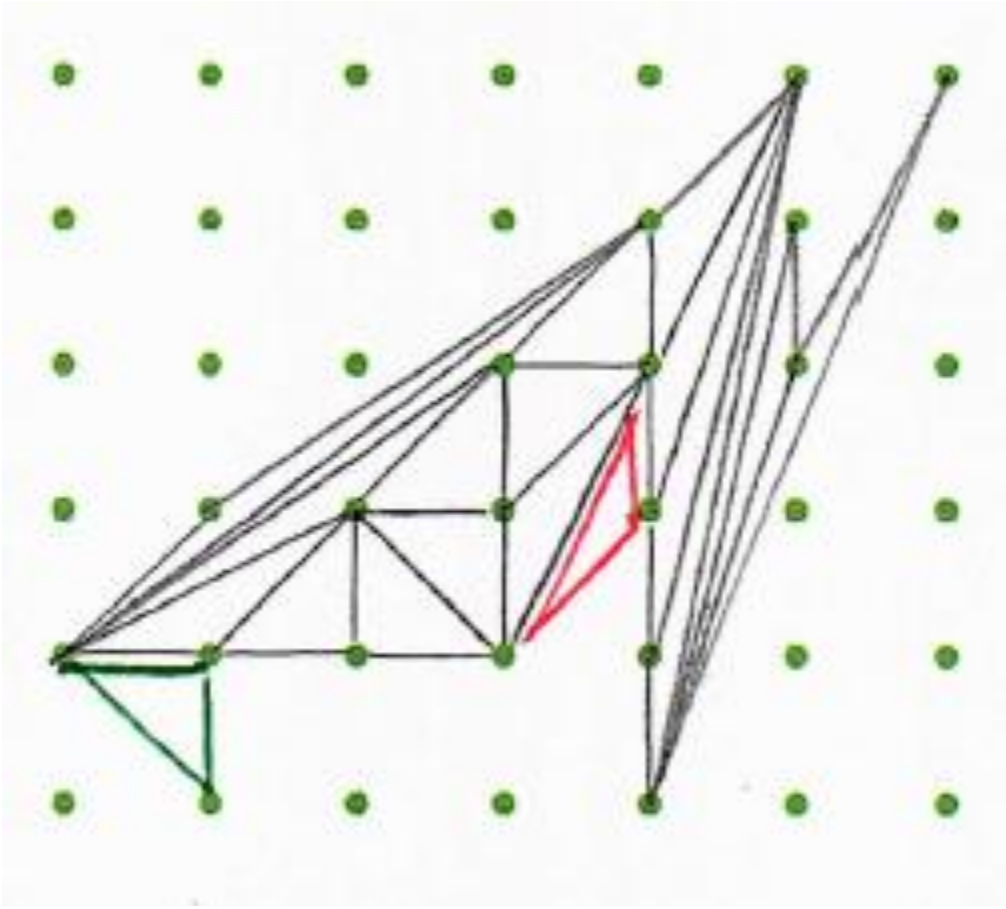

Fig,2

It means that our fully partitioned lattice polygon area can be re-assembled from these elementary triangles in a finite process validating Pick's formula. One also concludes, as a by-product, that the formula will still work for Jordan polygons assembled from triangles of unit area coming from still possible *aperiodic* tilings of plane by elementary triangles. To summarize, all we need in all approaches is to justify the following:

**Claim**: For every simple lattice vector there exist another simple lattice vector such that the area of a lattice parallelogram spanned on these two vectors is equal one.

Let **u** be an arbitrary simple vector which, in one of the equivalent bases, has coprime integer components $u_1, u_2$. From the well-known Bézout lemma (see, e.g. [21,22]) (or rather, from its special case for two coprime numbers) if follows that there exist two other coprime integers $v_1, v_2$ such that

$$u_1 v_1 + u_2 v_2 = 1 \qquad (6)$$

The integers $v_1, v_2$ can be easily found via Euclidean algorithm [21,22]-see also example 4.b in Section 4. .

Now, for definiteness, let's assume that we are working with the base $\{\mathbf{p},\mathbf{q}\}$. Then for any of two vectors $\pm(v_2\mathbf{p}-v_1\mathbf{q})$ parallelogram spanned on **u** and one of these vectors is a primitive cell of the lattice:

$$|(\pm(v_2\mathbf{p}-v_1\mathbf{q}))\times\mathbf{s}|=|(\pm(v_2\mathbf{p}-v_1\mathbf{q}))\times(u_1\mathbf{p}+u_2\mathbf{q})|= A(\Lambda).$$

This gives as an effective way of finding infinitely many different, equivalent bases (primitive cells) for a given two-dimensional lattice.

The example in Section 3.b illustrates how easily one can find a primitive cell of area 1 (or as cell declared as having a unit area for an arbitrary lattice) for given components of one of the two lattice simple vectors forming that cell. Calculations related to this task executed online, using nice application for Bézout relation [ 23].

## 3. On Pick's approach.

Pick's paper [1] entitled: *Geometrisches zur Zahlenlehre* (Geometric approach to the number theory) was published as a "re-worked' version of his lecture. The title and content is important, since Pick combines his measure having the additive property with a special triangulation argument, *to prove* an identity equivalent to the Bézout identity. He shows, *modo geometrico*, that for any two integers *a,b* with gcd(*a,b*)=*m* there exist integers $\alpha,\beta$ such that : $a\beta-b\alpha=m$ which is just the Bézout lemma in the "determinantial" form. He calls this result, with some exaggeration, *der Fundamental Satz der Zahlenteorie*. In our terminology, his proof proceeds from considering a lattice vector (*a,b*) passing through, by assumption, *m*-1 additional. "internal" lattice points. Next, a lattice triangle is formed having (*a,b*) as one side and some lattice vector $(\alpha_0,\beta_0)$ as its other side. This triangle has *at least* *m*+2 lattice points on the boundary and, typically, some lattice points in its interior. Its area is given by $\frac{1}{2}|a\beta_0-b\alpha_0|$ and, from Pick's formula, cannot be smaller than $\frac{1}{2}m$. This lower bound must be attained for some "internal", "smaller "triangle which does not contain any internal lattice points.

## 4. Examples, related problems, and comments.

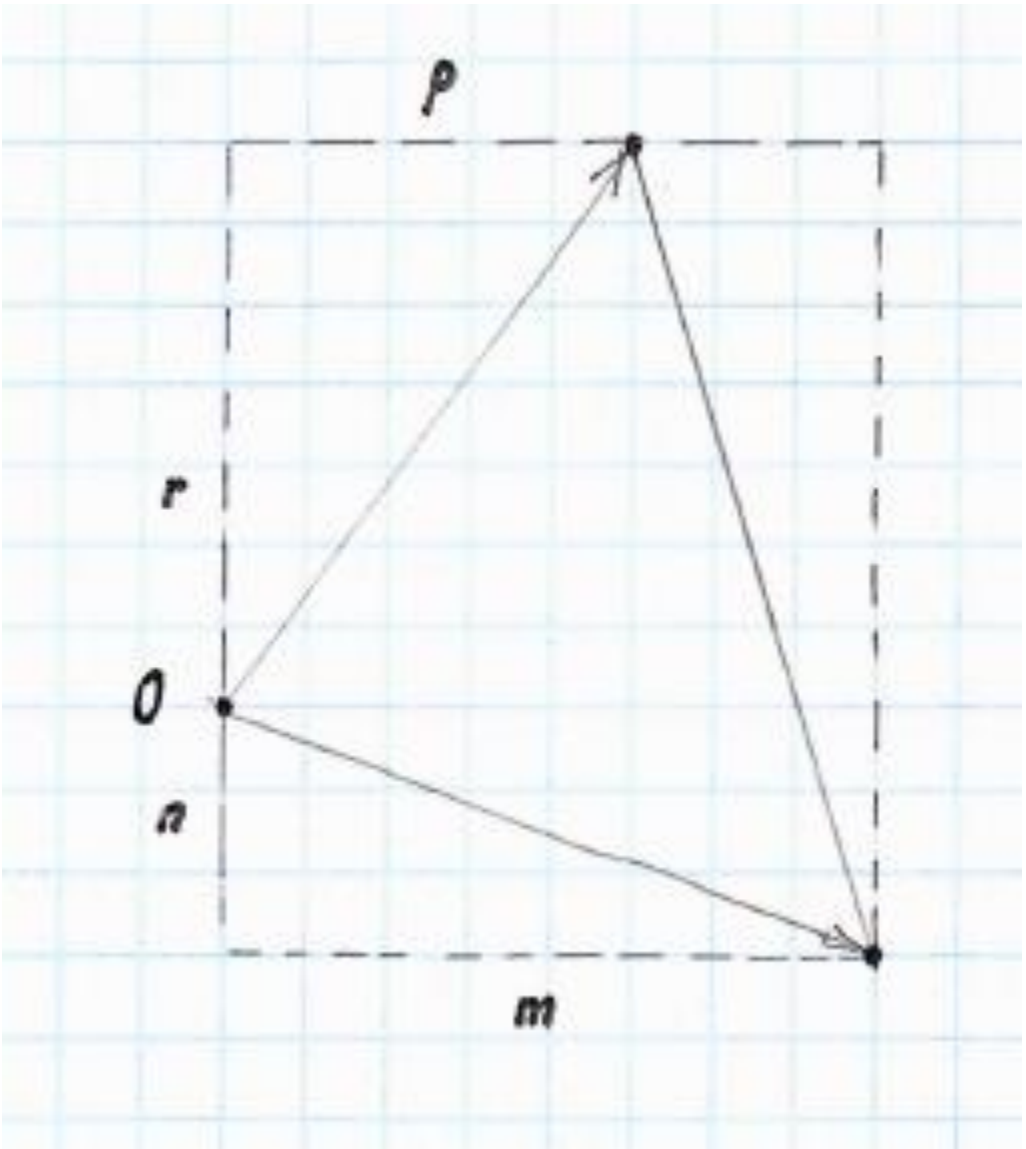

4.a. Lattice triangles (an elementary exercise)

Fig.3

A lattice triangle spanned on two lattice vectors with components (*p*, *r*) and (*m*, -*n*). Its area can be obtained from Pick's formula, since one knows that the formula works for the rectangle, three "corner" right triangles and then uses the additivity property. On the other hand, a student unfamiliar with vector products may check by elementary calculations that the triangle area is indeed given by the formula:

$$\frac{1}{2}\begin{vmatrix} m & -n \\ p & r \end{vmatrix} = \frac{1}{2}(mr + mn).$$

4.b. Using Euclid algorithm – an example: Consider simple lattice vector with coprime components (173, 16). The Euclid algorithm in this case reduces to 173=10$\cdot$16+13, 16= 1$\cdot$13 +3, 13= 4$\cdot$3 +1. Inverting these calculations in a standard way : 1= 13- 4$\cdot$3= 13- 4(16-1$\cdot$13)= (173-10$\cdot$16)-4(16-173 + 10$\cdot$16) , one obtains 1=5$\cdot$173-54$\cdot$16. Second simple lattice vector forming a primitive cell can thus be selected as a vector with components (54, 5). Clearly, other choices like (-54, -5) , (173-54, 16-5)=(119, 11) and (-119, -11) will work as well. Respective primitive cells are already difficult to draw on a "square grid paper", but still there are infinitely many other primitive cells with arbitrarily large "major" diagonals (compare it again to observation in [2]).

4.c. Proofs of Pick's formula by Euler's polyhedral formula (for planar connected graphs) and a formula for a fully triangulated Jordan polygon. These are proofs based upon an assumption that a given polygon can indeed be triangulated into triangles of identical area (set as equal ½). For a fully triangulated polygon containing $N^b$ points on the boundary and $N^i$ triangulation points inside one can easily find (by a finite induction, or rather a re-assembling process as in Section 2 )) a formula for the number $E$ od distinct edges in the triangulation [8]:

$$E = 2N^b + 3N^i - 3.$$

The number of edges $E$, on the other hand, appears in the well-known Euler's formula: $V + F - E = 2$ in standard notation (numbers of vertices, number of "faces", number of edges) applied for a connected planar graph which.in this case is a "completely" triangulated polygon, with $V = N^b + N^i$, $F = N(\triangle) + 1$ , total number of triangles increased by 1. Eliminating $E$ from two relations one obtains: $N(\triangle) = N^b + 2N^i - 2$ which is Pick's formula considered as a count of the number of elementary triangles. The method appears less direct than presented in Section 2, however relations with the Euler's polyhedral formula are quite important for generalizations of the Pick's result (see below, 4k).

4d. Another additive function for lattice Jordan polygons.

Let's cover every lattice point *i* of a given lattice Jordan polygon by sufficiently small disk $D_i$ and consider ratios

$$\alpha_i = \frac{A(P \cap D_i)}{A(D_i)}$$

where $A(.)$ is the area of the region considered. It is clear that $\alpha_i$ is always equal one for the internal lattice points. If *i* is a boundary point then $\alpha_i$ is the angular measure (in units equal $2\pi$) of the circular segment shared by *P* and $D_i$ which can be described as the "normalized visibility angle" of the interior of *P* for an observer located at the point *i*. Expression

$$\sum_{i \in P} \alpha_i$$

is a polygon characteristics which evidently has the additivity property for partitions of lattice polygons (like *P* from Fig. 1 split into two, where boundary visibility angles for $P_1$ and $P_2$ add up along the shared path). Noticing that $\sum_{i \in P} \alpha_i$ is equal to $\sum_{i \in boundary\ of\ P} \alpha_i + N^i(P)$ , and using an elementary fact that our Jordan lattice polygon is just an $N^b$-gon with $\sum_{i \in boundary\ of\ P} \alpha_i =$ $\frac{(N^b(P)-2)\pi}{2\pi} = \frac{1}{2} N_b(P) - 1$ we see that $\sum_{i \in P} \alpha_i$ is nothing but the Pick's formula "under disguise". Identification with the polygon area was carried out in [14], in a standard way (rectangles, right triangles, triangles.) in [15], directly for $\sum_{i \in boundary\ of\ P} \alpha_i + N^i(P)$. and then appealing to triangulation argument. Measures like $\sum_{i \in P} \alpha_i$ were introduced earlier in studies of d-dimensional polytopes with vertices covered by d-dimensional balls of sufficiently small radius- see Section 4.i and quoted there monograph [30] with complete history of "lattice points enumerators".

4.e. Larger parallelograms spanned by simple lattice vectors.

Lattice parallelograms spanned on two simple lattice vectors always have $N^b$=4. For them Pick's formula reduces to an elegant expression $A(P) = N^i(P) + 1$. For the unions of *k* such side-sharing congruent simple parallelograms it generalizes to $N^i(\bigcup P) + k$ .

4.f. There are several papers presenting Pick's formula as an interesting subject to be discussed at high school or college level: [24-27]. Examples of very elementary exercises related to Pick's formula (at the middle, or even elementary school level) can be found online: http://nrich.maths.org/1867/index and https://nrich.maths.org/5441.

4.h. Scaling argument.

Together with a given lattice polygon *P* let us consider all scaled versions of that polygon where all side lengths are scaled by an integer factor *k*= 2,3,…. . The areas of such scaled polygons $P'$ are equal $A(P')=k^2A(P)$ and $N^b(P')=kN^b(P)$. From Pick's formula it follows that,

$\lim_{k\to\infty}\frac{N^i(P')}{k^2}=A(P)$ in agreement with intuitive expectation that one will obtain ever increasing estimates (from below!) of the area *A*(*P*) by refining grid of the original lattice and counting only the number of internal points of this rescaled grid within polygon *P*.

4.i. Pick's formula and the surveyor's area formula.

Any lattice Jordan polygon is fully determined by an ordered sequence of lattice points on its boundary. This leads to a conclusion that such sequence should also contain information about polygon's area and that it should be encoded in an additive functional erasing information about internal points created in forming unions of two polygons which only share a part of their boundaries (for illustration see again Fig. 1).

Knowing that Pick's formula can be based upon triangulation of a lattice Jordan polygon into elementary triangles one can use this observation to relate it to other result for the area of a Jordan polygon – the surveyor's formula. Select an arbitrary lattice point as the origin of coordinate system and form a sequence of lattice points visited on the boundary of a lattice polygon starting from an arbitrary point ($m_0$, $n_0$) and moving in the counterclockwise direction: $(m_1,n_1),...(m_i,n_i),(m_{i+1},n_{i+1}),...(m_{N^b},n_{N^b})$. By assumption, all vectors

$$(m_2-m_1,n_2-n_1),..,(m_{i+1}-m_i,n_{i+1}-n_i),...(m_1-m_{N^b},n_1-n_{N^b})$$

are simple and each of them is a side of an elementary triangle in a given triangulation. Let $(u_i,v_i)$ be an internal polygon point for which the triangle with vertices $(m_i,n_i),(m_{i+1},n_{i+1}).(u_i,v_i)$ is an elementary one, in a given triangulation into elementary triangles. Unit area of this triangle can still be written, quite generally, as

$$\frac{1}{2}\begin{vmatrix} m_{i+1}-u_i & n_{i+1}-v_i \\ m_i-u_i & n_i-v_i \end{vmatrix}.$$

Yet another expression for that area can be obtained by using two different expansions of an auxiliary 3x3 determinant (similarly as it was done in [23]):

$$\frac{1}{2}\begin{vmatrix}1 & u_i & v_i\\ 1 & m_i & n_i\\ 1 & m_{i+1} & n_{i+1}\end{vmatrix} = \frac{1}{2}\begin{vmatrix}0 & u_i - m_i & v_i - n_i\\ 0 & m_i - m_{i+1} & n_i - n_{i+1}\\ 1 & m_{i+1} & n_{i+1}\end{vmatrix} = \frac{1}{2}\{\begin{vmatrix}u_i & v_i\\ m_i & m_i\end{vmatrix} + \begin{vmatrix}m_i & n_i\\ m_{i+1} & n_{i+1}\end{vmatrix} + \begin{vmatrix}m_{i+1} & n_{i+1}\\ u_i & v_i\end{vmatrix}\}$$

(A student not familiar with determinants may check the identity via direct calculation…).

The three determinants in the last expression can be considered as contribution from the oriented triangle sides with one of them determined by points on the boundary. Now let's consider a complete triangulation of the polygon into *oriented* (in the assumed counterclockwise direction) elementary triangles with unit area given by the expression of the last type and sum all of them. Every internal (non-boundary) side of these triangles is then traversed twice in the *opposite directions* and contributions from these sides disappear in the sum, in yet another manifestation of the additivity property. The final result is the surveyor formula:

$$A(P) = \tfrac{1}{2} \sum_{i \in \text{ boundary of } P} \begin{vmatrix} m_i & n_i \\ m_{i+1} & n_{i+1} \end{vmatrix} ., \; i_{N^b+1} \equiv i_1 .$$

If three or more consecutive lattice points on the boundary of a lattice Jordan polygon are aligned in the same direction these "non-vertex" points can be eliminated from this expression (another easy exercise) leading to surveyor formula in standard form where summation goes over polygon vertices.

All of the above can be used as elementary, introductory or complementary material in calculus or mathematical methods classes when working with line integrals in real and complex planes and in a discussion of the Green' theorem.

4.j. The Farey sequence

This family of finite sequences (commonly called Farey series) was also a subject of $ 3 of Pick's paper where he used again, a geometric, lattice based approach to prove their basic properties. For given natural $n$=1,2,.. let's construct a finite sequence of all proper, non-negative fractions with denominators not exceeding $n$ and arranged in the ascending order. One basic property of this sequences is that for any two subsequent elements in such sequences, say $a/b < c/d$, one has $bc - ad = 1$. It follows from our previous discussion that elements of Farey sequences (with clearly co-prime numerator and denominator) can be used in construction of primitive cells and elementary triangles. The convoluted history of Farey sequences and their relationship with Pick's theorem was presented in [20]. This is an enjoyable review recommended to anyone interested in true "history of mathematical discovery".

4.k. Extensions and generalization of the Pick's theorem.

Relatively simple generalization of Pick's theorem to "lattice polygons with holes", special unions of lattice polygons, and non-standard lattices etc. (e.g., [15], [28].[34]). are scattered in cited below literature. Progress in this direction (until 1993) is reviewed in [31] and one important generalization added for polygons with overlapping edges.

Very interesting and studied for a long time are results related to generalizations of Pick's theorem for more challenging case of lattice polytopes in higher dimensions. Early results in this direction were obtained by Reeve [29] and progress in this direction until 1989 is summarized in an excellent monograph "*Lattice Points"* by Erdös and co-authors [30]. Pick's result was also obtained by using Minkowski's theorem in [32]. Quite recent results on sums of lattice points were also obtained in [33].

Finally, fascinating connections with problems of modern theoretical physics were demonstrated in [35] -"but that's another story"...

*) DEPARTMENT OF PHYSICS, UNIVERSITY OF NORTH TEXAS
-retired Assoc. Prof. of Physics
*E-mail address:*jacek.kowalski@unt.edu